\documentclass[11pt]{amsart}
\usepackage{amsmath,amssymb,amsthm}
\usepackage[margin=1.15in]{geometry}
\usepackage[colorlinks,linkcolor=blue,citecolor=blue]{hyperref}

\newtheorem{thm}{Theorem}[section]
\newtheorem{lema}[thm]{Lemma}
\newtheorem{prop}[thm]{Proposition}

\theoremstyle{definition}
\newtheorem{defn}[thm]{Definition}
\newtheorem{exam}[thm]{Example}
\theoremstyle{remark}
\newtheorem{rem}[thm]{Remark}
\numberwithin{equation}{section}

\def\Xint#1{\mathchoice
  {\XXint\displaystyle\textstyle{#1}}%
  {\XXint\textstyle\scriptstyle{#1}}%
  {\XXint\scriptstyle\scriptscriptstyle{#1}}%
  {\XXint\scriptscriptstyle\scriptscriptstyle{#1}}%
  \!\int}
\def\XXint#1#2#3{{\setbox0=\hbox{$#1{#2#3}{\int}$}\vcenter{\hbox{$#2#3$}}\kern-.5\wd0}}
\newcommand{\fint}{\Xint-}

\begin{document}

\title{Fractional order Orlicz--Sobolev spaces revisited}

\author[J. Fern\'andez Bonder]{Juli\'an Fern\'andez Bonder}
\address{Departamento de Matem\'atica, Facultad de Ciencias Exactas y Naturales, Universidad de Buenos Aires, Buenos Aires, Argentina, and Instituto de C\'alculo UBA-CONICET, Buenos Aires, Argentina}
\email{jfbonder@dm.uba.ar}

\author[A. M. Salort]{Ariel M. Salort}
\address{Departamento de Matem\'aticas Aplicadas y Ciencia de Datos, Universidad San Pablo-CEU, CEU Universities, Juli\'an Romea 2, 28003 Madrid, Spain}
\email{ariel.salort@ceu.es}

\subjclass[2020]{46E30, 35R11, 49J45, 35J92}
\keywords{Fractional order Orlicz--Sobolev spaces, $\Delta_2$ condition, Bourgain--Brezis--Mironescu formula, $\Gamma$-convergence, $g$-Laplacian}

\begin{abstract}
The aim of this paper is to revisit our work [J. Funct. Anal. 277 (2019), 333--367], where fractional order Orlicz--Sobolev spaces were introduced and their behavior as the fractional parameter $s\uparrow1$ was studied, and to redo it without assuming the $\Delta_2$ condition on the Young function $G$ or on its complementary function. Since that paper appeared, the theory of these spaces has been considerably developed, and the tools available today allow for a cleaner treatment in full generality. We prove a Bourgain--Brezis--Mironescu type formula as $s\uparrow1$, a compactness and lower semicontinuity result for sequences with bounded energy, a Poincar\'e inequality uniform in $s$, and the $\Gamma$-convergence of the fractional modulars to a local modular on arbitrary open sets. As an application we study the asymptotic behavior of minimizers and of solutions of Dirichlet problems for the fractional $g$-Laplacian as $s\uparrow1$. This version also corrects some errors in the original paper, and it is intended to serve as a more complete and reliable reference for these results.
\end{abstract}

\maketitle

\section{Introduction}

Fractional order Sobolev spaces $W^{s,p}$ are the natural functional framework for many nonlocal problems, and a celebrated result of Bourgain, Brezis and Mironescu \cite{BBM} describes their behavior as $s\uparrow1$: for $u\in L^p(\mathbb{R}^n)$,
\[
\lim_{s\uparrow1}(1-s)\iint_{\mathbb{R}^n\times\mathbb{R}^n}\frac{|u(x)-u(y)|^p}{|x-y|^{n+sp}}\,dx\,dy=\frac{K_{n,p}}p\int_{\mathbb{R}^n}|\nabla u|^p\,dx ,
\]
where $K_{n,p}=\int_{\mathbb S^{n-1}}|\omega\cdot e_n|^p\,d\mathcal H^{n-1}$ (see Section \ref{sec.prelim}).
When the power $a^p$ is replaced by a general Young function $G$, one is led to the fractional Orlicz--Sobolev spaces introduced in \cite{FBS}, defined through the modular
\[
\iint_{\mathbb{R}^n\times\mathbb{R}^n}G\Big(\frac{|u(x)-u(y)|}{|x-y|^s}\Big)\frac{dx\,dy}{|x-y|^n}.
\]
In \cite{FBS} the asymptotic behavior of these spaces as $s\uparrow1$ was studied under the $\Delta_2$ condition $G(2a)\le\mathbf CG(a)$. This condition excludes Young functions with fast growth, such as $G(a)=e^{a^2}-1$ or $G(a)=e^a-1-a$, and it is used in \cite{FBS} in several essential ways: the Orlicz class is a vector space and coincides with $L_G$, smooth functions are dense in norm, the spaces are separable, and modular and norm convergence are equivalent.

The purpose of this paper is to redo \cite{FBS} without the $\Delta_2$ condition. Several years have passed since \cite{FBS} was written, and during this time the theory of fractional Orlicz--Sobolev spaces has been developed by many authors, both with and without growth conditions on $G$. The tools available today allow for a treatment of the results of \cite{FBS} in full generality and with simpler proofs, and we believe it is worthwhile to collect them in a single place. Some results in this direction have been obtained recently: Alberico, Cianchi, Pick and Slav\'{\i}kov\'a \cite{ACPS} studied the limit as $s\uparrow1$ of possibly non-separable fractional Orlicz--Sobolev spaces; the $\Gamma$-convergence of the modulars without $\Delta_2$ was obtained in \cite{CFBS}; the asymptotic behavior as $s\uparrow1$ of fractional Musielak--Orlicz--Sobolev modulars without the $\Delta_2$ condition was studied in \cite{AACS}; and eigenvalue problems for the fractional $g$-Laplacian without $\Delta_2$ were treated in \cite{FBSp,SV}, see also \cite{FBGS} for nonlocal anisotropic operators.  More generally, fractional Orlicz--Sobolev spaces built upon arbitrary Young functions have been systematically studied by Alberico, Cianchi, Pick and Slav\'{\i}kov\'a, who obtained optimal Sobolev embeddings and compactness criteria \cite{ACPS2}, necessary and sufficient conditions for the boundedness and continuity of their elements \cite{ACPS4}, and the limit as $s\downarrow0$ \cite{ACPS0}, where a Maz'ya--Shaposhnikova type formula is proved under the $\Delta_2$ condition, which is shown there to be essential; see also the survey \cite{ACPS3}. For the general theory of Orlicz and Musielak--Orlicz spaces without growth conditions we refer to \cite{HH, CGSW}, and for the approximation of Orlicz--Sobolev functions in the absence of $\Delta_2$ to \cite{Gossez, Gossez82}.  Here we give a unified and self-contained treatment of the results of \cite{FBS} for arbitrary Young functions, building on these developments.

This paper also corrects some errors in \cite{FBS}, and it is intended to replace it as a reference for these results. Lemma 5.3 of \cite{FBS} does not hold as stated (see Remark \ref{rem.lema53}); it was used in the proof of the uniform estimate \cite[Theorem 5.2]{FBS}, which is replaced here by Theorem \ref{teo.s1.s2}. Moreover, the existence of the limit defining the function $G_0$ was assumed in \cite{FBS}, while it always exists (see \eqref{exacto}). All the results of \cite{FBS} remain valid, and in fact they hold in greater generality.

Our main results are the following:
\begin{itemize}
\item the limit Young function $G_0$ always exists and is equivalent to $G$ in the sense of Orlicz spaces (Proposition \ref{prop.equiv}), although in general not with multiplicative constants;
\item the compact embedding of $W^sL_G(\mathbb{R}^n)$ into $E_{G,{\rm loc}}(\mathbb{R}^n)$ (Theorem \ref{teo.comp});
\item the Bourgain--Brezis--Mironescu formula (Theorem \ref{main1});
\item a uniform estimate comparing the modulars of orders $t\le s$ (Theorem \ref{teo.s1.s2}), and the compactness and lower semicontinuity of sequences with bounded energy as $s\uparrow1$ (Theorems \ref{main2} and \ref{teo.liminf});
\item a Poincar\'e inequality uniform in $s$ (Theorem \ref{teo.poinc.unif}) and the $\Gamma$-convergence of the fractional modulars on arbitrary open sets (Theorem \ref{teo.gamma.conv});
\item the convergence of minimizers and of (variational) solutions of Dirichlet problems for the fractional $g$-Laplacian to the solution of the corresponding local problem (Theorems \ref{teo.conv.min} and \ref{teo.conv.sol}).
\end{itemize}
A structural feature of the present approach is the limit Young function $G_0$. It is not merely the integrand appearing in the local limit: Proposition \ref{prop.equiv} shows that it is always a Young function and is equivalent to $G$ after a fixed rescaling of the argument, while Lemma \ref{lema.Gt.reg} shows that $G_0\in C^1([0,\infty))$ and is strictly convex even when $G$ has neither property. This regularization is important in the local variational problem, where it yields differentiability and strict convexity without imposing (G1) or (G2) on $G$ (see Section \ref{sec.glap}).
No condition besides being a Young function is imposed on $G$, except in Section \ref{sec.glap}, where the regularity of $G$ is assumed to derive the Euler--Lagrange equation of the fractional problem, and its strict convexity to have uniqueness of its solutions; the limit problem requires neither assumption. In particular, no growth condition, neither from above ($\Delta_2$) nor from below ($\nabla_2$), is assumed. Moreover, no regularity is assumed on the open sets involved.
Since the growth assumptions are absent, it is essential throughout to distinguish the Orlicz class $\mathcal L_G$, the Orlicz space $L_G$, and its separable subspace $E_G$. The same care is needed for the Sobolev spaces: in particular, $W^s_0L_G(\Omega)$ is defined by zero extension outside $\Omega$, whereas $W^s_0E_G(\Omega)=W^s_0L_G(\Omega)\cap W^sE_G(\mathbb{R}^n)$. The effective domain of the limit functional is the modular domain described in Remark \ref{rem.ae}, and should not be identified in general with $W^1_0E_G(\Omega)$.

The paper is organized as follows. In Section \ref{sec.prelim} we collect the necessary preliminaries on Orlicz spaces and prove some technical lemmas. Section \ref{sec.comp} is devoted to the compactness theorem. In Section \ref{sec.bbm} we prove the Bourgain--Brezis--Mironescu formula, and in Section \ref{sec.seq} we treat sequences. Finally, Section \ref{sec.app} contains the applications.

\section{Preliminaries}\label{sec.prelim}

\subsection{Young functions}

\begin{defn}
A function $G\colon[0,\infty)\to[0,\infty)$ is a \emph{Young function} if
\[
G(a)=\int_0^ag(\tau)\,d\tau,
\]
where $g\colon[0,\infty)\to[0,\infty)$ is nondecreasing and right continuous, $g(0)=0$, $g(a)>0$ for $a>0$ and $g(a)\to\infty$ as $a\to\infty$.
\end{defn}

A Young function is convex, increasing and locally Lipschitz, $G(0)=0$, and
\[
\lim_{a\to0^+}\frac{G(a)}{a}=0,\qquad\lim_{a\to\infty}\frac{G(a)}{a}=\infty .
\]
Convexity gives
\begin{equation}\label{P5}
G(\tau a)\le\tau G(a)\quad(0\le\tau\le1),\qquad G(\tau a)\ge\tau G(a)\quad(\tau\ge1),
\end{equation}
and
\begin{equation}\label{P4}
G(a+b)\le\tfrac12\big(G(2a)+G(2b)\big),\qquad a,b\ge0 .
\end{equation}
Note that \eqref{P4} replaces the inequality $G(a+b)\le\frac{\mathbf C}{2}(G(a)+G(b))$ used in \cite{FBS}, which requires the $\Delta_2$ condition.

The complementary function is $\overline{G}(b)=\sup\{ab-G(a)\colon a\ge0\}$. It is also a Young function, and Young's inequality $ab\le G(a)+\overline{G}(b)$ holds. Moreover, $\overline{G}(g(a))=ag(a)-G(a)$ and, since $ag(a)\le\int_a^{2a}g\le G(2a)$,
\begin{equation}
\overline{G}(g(a))\le G(2a),\qquad a\ge0 .
\end{equation}

We say that $G$ satisfies the $\Delta_2$ condition if $G(2a)\le\mathbf CG(a)$ for all $a\ge0$, and that it satisfies the $\Delta_2$ condition near infinity if this inequality holds for all $a\ge a_0$, for some $a_0>0$. \emph{Neither condition is assumed in this paper.}

\subsection{Orlicz spaces}
We follow \cite{KR, Gossez} and \cite{RR}; see also \cite{HH, CGSW}. Let $(U,\mu)$ be a measure space. For $u$ measurable we define the modular
\[
\Phi_{G,\mu}(u)=\int_UG(|u|)\,d\mu ,
\]
and we simply write $\Phi_G$ when $\mu$ is the Lebesgue measure on $\mathbb{R}^n$, that is, $\Phi_G(u)=\int_{\mathbb{R}^n} G(|u|)\,dx$. The \emph{Orlicz class} is
\[
\mathcal L_G(U,\mu)=\{u\colon \Phi_{G,\mu}(u)<\infty\},
\]
which, for the Lebesgue measure on $\mathbb{R}^n$, is a vector space if and only if $G$ satisfies $\Delta_2$. The \emph{Orlicz space} $L_G(U,\mu)$ is the linear hull of $\mathcal L_G(U,\mu)$, that is, the set of measurable $u$ with $\Phi_{G,\mu}(u/k)<\infty$ for some $k>0$. It is a Banach space with the Luxemburg norm
\[
\|u\|_{G,\mu}=\inf\Big\{k>0\colon\Phi_{G,\mu}\Big(\frac{u}{k}\Big)\le1\Big\},
\]
and we write $\|u\|_G$ when $\mu$ is clear from the context. Recall that, for $v\in L_G(U,\mu)$,
\begin{equation}\label{mod-norm}
\|v\|_G\le1\ \Rightarrow\ \Phi_{G,\mu}(v)\le\|v\|_G,\qquad\Phi_{G,\mu}(v)\le C\ \Rightarrow\ \|v\|_G\le\max\{1,C\}.
\end{equation}
The subspace $E_G(U,\mu)$ is the closure in $L_G$ of the bounded functions with support of finite measure, and
\[
E_G(U,\mu)=\{u\colon\Phi_{G,\mu}(u/k)<\infty\text{ for every }k>0\}.
\]
We have $E_G\subset\mathcal L_G\subset L_G$, with equalities if and only if $G$ satisfies $\Delta_2$ in the appropriate sense: globally when $\mu$ is the Lebesgue measure on $\mathbb{R}^n$ (or on an open set of infinite measure), and near infinity when $\mu$ is the Lebesgue measure on an open set $U\subset\mathbb{R}^n$ of finite measure. The space $E_G$ is separable, $(E_G)^*=L_{\overline{G}}$, and $L_G$ is reflexive if and only if both $G$ and $\overline{G}$ satisfy $\Delta_2$ (in the same sense). H\"older's inequality reads $\int|uv|\,d\mu\le2\|u\|_G\|v\|_{\overline{G}}$.

When $U\subset\mathbb{R}^n$ and $\mu$ is the Lebesgue measure, we write $L_G(U)$, $E_G(U)$.

\subsection{Orlicz--Sobolev spaces}
On $\mathbb{R}^{2n}$ we consider the measure
\[
d\nu=\frac{dx\,dy}{|x-y|^n}.
\]
For $s\in(0,1)$ and $u$ measurable we write $D^su(x,y)=\frac{u(x)-u(y)}{|x-y|^s}$, and we define the fractional modular
\[
\Phi_{G,\nu}(D^su)=\iint_{\mathbb{R}^{2n}}G(|D^su|)\,d\nu .
\]
The fractional Orlicz--Sobolev spaces are
\begin{align*}
W^sL_G(\mathbb{R}^n)&=\{u\in L_G(\mathbb{R}^n)\colon D^su\in L_G(\mathbb{R}^{2n},\nu)\},\\
W^sE_G(\mathbb{R}^n)&=\{u\in E_G(\mathbb{R}^n)\colon D^su\in E_G(\mathbb{R}^{2n},\nu)\}.
\end{align*}
In order to have the correct scaling as $s\uparrow1$, we endow them with the norms defined through the normalized measures
\[
\nu_s=(1-s)\,\nu,\qquad\text{for which}\qquad\Phi_{G,\nu_s}(D^su)=(1-s)\,\Phi_{G,\nu}(D^su).
\]
Namely, we let
\[
\|u\|_{s,G}=\|u\|_G+[u]_{s,G},\qquad[u]_{s,G}=\|D^su\|_{G,\nu_s}.
\]
Since $\Phi_{G,\nu_s}=(1-s)\Phi_{G,\nu}$, convexity gives
\begin{equation}
(1-s)\|D^su\|_{G,\nu}\le[u]_{s,G}\le\|D^su\|_{G,\nu},
\end{equation}
so for each fixed $s$ this norm is equivalent to the one defined with $\nu$, and the spaces do not change. The advantage of $[u]_{s,G}$ is that it remains bounded as $s\uparrow1$ for smooth $u$ (see Theorem \ref{main1} and Lemma \ref{lema.cota}).
The space $W^sL_G(\mathbb{R}^n)$ is a Banach space and $W^sE_G(\mathbb{R}^n)$ is a closed subspace. In the same way we define $W^1L_G(\mathbb{R}^n)$ and $W^1E_G(\mathbb{R}^n)$, requiring $u$ and $|\nabla u|$ to belong to $L_G$ (resp.\ $E_G$).

\paragraph{Functions vanishing outside a set.} Let $\Omega\subset\mathbb{R}^n$ be open. Functions defined on $\Omega$ are extended by zero, and for $s\in(0,1]$ we define
\[
W^s_0L_G(\Omega)=\{u\in W^sL_G(\mathbb{R}^n)\colon u=0\text{ a.e.\ in }\mathbb{R}^n\setminus\Omega\},
\]
and $W^s_0E_G(\Omega)=W^s_0L_G(\Omega)\cap W^sE_G(\mathbb{R}^n)$. No regularity of $\Omega$ is assumed. By Lemma \ref{lema.cota} below, $C^\infty_c(\Omega)\subset W^s_0E_G(\Omega)$.

\begin{rem}
In \cite{FBSp}, $W^s_0L_G(\Omega)$ is defined as the closure of $C^\infty_c(\Omega)$ in the weak-$*$ topology, which requires the segment property on $\Omega$ for several results. The definitions above do not require any regularity of $\Omega$ and, as shown in \cite{CFBS}, they are the natural ones for the $\Gamma$-convergence as $s\uparrow1$ (see Remark \ref{rem.ae}).
\end{rem}

\subsection{Regularization}
Let $\rho\in C^\infty_c(B_1)$ be a standard mollifier and $\rho_\varepsilon(x)=\varepsilon^{-n}\rho(x/\varepsilon)$. For $u\in L_G(\mathbb{R}^n)$ let $u_\varepsilon=u*\rho_\varepsilon$.

\begin{lema}\label{lema.C2}
Let $u\in L_G(\mathbb{R}^n)$ and $\varepsilon>0$. Then $u_\varepsilon\in C^\infty(\mathbb{R}^n)$ and, for every multi-index $\alpha$,
\[
\|\partial^\alpha u_\varepsilon\|_\infty\le2\|\partial^\alpha\rho_\varepsilon\|_{\overline{G}}\|u\|_G .
\]
Moreover, if $u_k\to u$ in $L^1_{\rm loc}(\mathbb{R}^n)$, then $(u_k)_\varepsilon\to u_\varepsilon$ in $C^1(K)$ for every compact $K$.
\end{lema}

\begin{proof}
The bound follows from $\partial^\alpha u_\varepsilon(x)=\int_{\mathbb{R}^n} u(y)\,\partial^\alpha\rho_\varepsilon(x-y)\,dy$, H\"older's inequality and the translation invariance of $\|\cdot\|_{\overline{G}}$. The second assertion follows from the same formula, since $\partial^\alpha\rho_\varepsilon$ is bounded and supported in $B_\varepsilon$.
\end{proof}

\begin{lema}\label{lema.reg}
Let $u\in L_G(\mathbb{R}^n)$. Then $\Phi_{G,\nu_s}(D^su_\varepsilon)\le\Phi_{G,\nu_s}(D^su)$ for every $\varepsilon>0$ and $0<s<1$.
\end{lema}

\begin{proof}
This is \cite[Lemma 2.13]{FBS}; the proof uses only Jensen's inequality and Fubini's theorem.
\end{proof}

\begin{lema}\label{lema.translation}
Let $0<s<1$ and $u\in W^sL_G(\mathbb{R}^n)$. Then, for every $h\in\mathbb{R}^n$ with $|h|\le1$,
\[
\|\tau_hu-u\|_G\le C_n\,|h|^s\,[u]_{s,G},\qquad C_n=\frac{32\cdot3^n}{\omega_n},
\]
where $\tau_hu(x)=u(x+h)$ and $\omega_n=|B_1|$. The constant does not depend on $s$.
\end{lema}

\begin{proof}
By homogeneity we assume $[u]_{s,G}=1$, that is, $\Phi_{G,\nu}(D^su)\le\frac1{1-s}$. Let $r_0=|h|\in(0,1]$.

\emph{Step 1: control of the modular at scale $r_0$.} For $r>0$ and $\varrho\ge0$ let
\[
\Psi(r,\varrho)=\int_{\mathbb S^{n-1}}\int_{\mathbb{R}^n} G\Big(\varrho\,\frac{|u(x+r\omega)-u(x)|}{r}\Big)dx\,d\mathcal H^{n-1}(\omega).
\]
Since $\frac{|u(x+2r\omega)-u(x)|}{2r}$ is bounded by the average of the corresponding quotients at scale $r$ at the points $x$ and $x+r\omega$, convexity and translation invariance give $\Psi(2r,\varrho)\le\Psi(r,\varrho)$. Let $I_j=[2^{-j-1}r_0,2^{-j}r_0]$ and $\varphi_j(\varrho)=\int_{I_j}\Psi(r,\varrho)\frac{dr}{r}$. The change of variables $r=2^{-j}r'$ gives $\varphi_j\ge\varphi_0$ for every $j\ge0$. Let $M=\lceil\frac1{1-s}\rceil$. For $0\le j\le M-1$ we have $(j+1)(1-s)\le3$, and since $r^{1-s}\ge2^{-(j+1)(1-s)}r_0^{1-s}\ge\frac18r_0^{1-s}$ on $I_j$,
\[
M\varphi_0\Big(\frac{r_0^{1-s}}8\Big)\le\sum_{j=0}^{M-1}\varphi_j\Big(\frac{r_0^{1-s}}8\Big)\le\int_0^{r_0}\Psi(r,r^{1-s})\frac{dr}r\le\Phi_{G,\nu}(D^su)\le\frac1{1-s}\le M .
\]
Since $r_0^{1-s}/r\ge r^{-s}$ for $r\le r_0$, we conclude that
\[
\int_{r_0/2}^{r_0}\int_{\mathbb S^{n-1}}\int_{\mathbb{R}^n} G\Big(\frac{|u(x+r\omega)-u(x)|}{8r^s}\Big)dx\,d\mathcal H^{n-1}\,\frac{dr}r\le\varphi_0\Big(\frac{r_0^{1-s}}8\Big)\le1 .
\]
The same holds with $r_0$ replaced by $r_0/2$, so the integral over $r\in[r_0/4,r_0]$ is at most $2$.

\emph{Step 2: averaging.} Let $B=B(h/2,r_0/4)$. For $z\in B$, both $|z|$ and $|h-z|$ belong to $[r_0/4,3r_0/4]$, and $z\mapsto h-z$ maps $B$ onto itself. Writing $u(x+h)-u(x)=(u(x+h)-u(x+z))+(u(x+z)-u(x))$, averaging over $z\in B$ and using the convexity of $G$ and translation invariance, for every $\mu>0$,
\[
\int_{\mathbb{R}^n} G\Big(\frac{|u(x+h)-u(x)|}{2\mu}\Big)dx\le\fint_BH_\mu(w)\,dw,\qquad H_\mu(w)=\int_{\mathbb{R}^n} G\Big(\frac{|u(x+w)-u(x)|}{\mu}\Big)dx .
\]
In polar coordinates, since $r^n\le(3r_0/4)^n$ for $r\le3r_0/4$,
\[
\fint_BH_\mu\,dw\le\frac{3^n}{\omega_n}\int_{r_0/4}^{r_0}\int_{\mathbb S^{n-1}}H_\mu(r\omega)\,d\mathcal H^{n-1}\,\frac{dr}r .
\]
Take $\mu=8Kr_0^s$ with $K=\frac{2\cdot3^n}{\omega_n}\ge1$. For $r\le r_0$, \eqref{P5} gives $G\big(\frac{a}{8Kr_0^s}\big)\le\frac1KG\big(\frac a{8r^s}\big)$, so by Step 1 the right-hand side is at most $\frac{3^n}{\omega_n}\cdot\frac2K=1$. Hence $\|\tau_hu-u\|_G\le2\mu=16Kr_0^s=C_n|h|^s$.
\end{proof}

\begin{rem}
For $s=1$, Jensen's inequality gives $\|\tau_hu-u\|_G\le|h|\,\|\nabla u\|_G$. Thus Lemma \ref{lema.translation}, whose constant does not depend on $s$, is consistent with the limit case. Compare with \cite[Lemma B.1, Corollary B.2]{FBSp}, where a similar estimate is proved with constant depending on $s$.
\end{rem}

\subsection{The limit Young function}

For $a\ge0$ let
\begin{equation}
G_0(a)=\int_0^1\int_{\mathbb S^{n-1}}G(a\,r\,|\omega\cdot e_n|)\,d\mathcal H^{n-1}(\omega)\,\frac{dr}{r}=\int_0^a\int_{\mathbb S^{n-1}}G(r|\omega\cdot e_n|)\,d\mathcal H^{n-1}(\omega)\,\frac{dr}{r}.
\end{equation}
In \cite{FBS}, this function (denoted $\tilde G$ there) was defined as the limit (when it exists) of $(1-s)\int_0^1\int_{\mathbb S^{n-1}}G(a|\omega\cdot e_n|r^{1-s})\,d\mathcal H^{n-1}\,\frac{dr}{r}$ as $s\uparrow1$. The change of variables $\tau=r^{1-s}$, for which $(1-s)\frac{dr}{r}=\frac{d\tau}{\tau}$, shows that this quantity does not depend on $s$:
\begin{equation}\label{exacto}
(1-s)\int_0^1\int_{\mathbb S^{n-1}}G(a|\omega\cdot e_n|r^{1-s})\,d\mathcal H^{n-1}\,\frac{dr}{r}=G_0(a)\qquad\text{for every }s\in(0,1).
\end{equation}
In particular, the limit always exists, and the assumption of its existence in \cite{FBS} can be removed.

\begin{prop}\label{prop.equiv}
$G_0$ is a Young function, $\xi\mapsto G_0(|\xi|)$ is convex in $\mathbb{R}^n$, and
\[
c_n\,G(a/4)\le G_0(a)\le n\omega_n\,G(a)\qquad\text{for all }a\ge0,
\]
where $c_n=\log2\cdot\mathcal H^{n-1}(\{\omega\in\mathbb S^{n-1}\colon|\omega\cdot e_n|\ge\frac12\})$. Consequently, there are constants $0<\alpha_n\le1\le\beta_n$ such that
\begin{equation}\label{equiv.escala}
G(\alpha_na)\le G_0 (a)\le G(\beta_na)\qquad\text{for all }a\ge0,
\end{equation}
and $L_{ G_0 }=L_G$, $E_{ G_0 }=E_G$ with equivalent norms.
\end{prop}

\begin{proof}
Let $\Lambda(\tau)=\int_{\mathbb S^{n-1}}G(\tau|\omega\cdot e_n|)\,d\mathcal H^{n-1}$, which is convex, nondecreasing and satisfies $\Lambda(\tau r)\le r\Lambda(\tau)$ for $r\in[0,1]$. Then $ G_0 (a)=\int_0^1\Lambda(ar)\frac{dr}{r}$ is finite, convex and nondecreasing, so $ G_0 (|\xi|)$ is convex. Since $\Lambda(r)\le n\omega_nG(r)$ and $\int_0^1G(ar)\frac{dr}{r}\le G(a)$ by \eqref{P5}, the upper bound follows. Restricting the integral to $r\in[\frac12,1]$ and $|\omega\cdot e_n|\ge\frac12$ gives the lower bound. Writing $ G_0 (a)=\int_0^a\frac{\Lambda(r)}{r}dr$, its derivative $\Lambda(a)/a$ is continuous and nondecreasing (since $\Lambda$ is convex with $\Lambda(0)=0$), it tends to $0$ as $a\to0$ (since $\Lambda(a)\le n\omega_nG(a)$) and to $\infty$ as $a\to\infty$ (since $\Lambda(a)\ge\mathcal H^{n-1}(\{|\omega\cdot e_n|\ge\frac12\})\,G(a/2)$). Hence $ G_0 $ is a Young function.

Finally, \eqref{equiv.escala} follows from \eqref{P5} with $\alpha_n=\frac14\min\{1,c_n\}$ and $\beta_n=\max\{1,n\omega_n\}$: indeed $G(\alpha_na)\le\min\{1,c_n\}G(\frac a4)\le c_nG(\frac a4)\le G_0 (a)$ and $ G_0 (a)\le n\omega_nG(a)\le\beta_nG(a)\le G(\beta_na)$.
\end{proof}

On $W^1L_G(\mathbb{R}^n)$ we use the norm
\[
\|u\|_{1,G}=\|u\|_G+[u]_{1,G},\qquad[u]_{1,G}=\|\nabla u\|_G,
\]
where $\|\nabla u\|_G$ denotes the Luxemburg norm of $|\nabla u|$ in $L_G(\mathbb{R}^n)$. Analogously $[u]_{1, G_0 }=\|\nabla u\|_{ G_0 }$. Since $\Phi_{G,\nu_s}(D^su)\to\Phi_{ G_0 }(|\nabla u|)$ as $s\uparrow1$ (Theorem \ref{main1}), the natural limit of $[u]_{s,G}$ is $[u]_{1, G_0 }$. By \eqref{equiv.escala}, we have that following:
\begin{equation}\label{equiv.normas.1}
\alpha_n[u]_{1,G}\le[u]_{1, G_0 }\le\beta_n[u]_{1,G}.
\end{equation}

\begin{rem}
Without the $\Delta_2$ condition, $G$ and $ G_0 $ are \emph{not} equivalent with multiplicative constants in front, that is, $G\le C G_0 $ may fail. For instance, let $G(a)=e^{a^2}-1$. Since $\Lambda(r)\le n\omega_nG(r)$, we have $ G_0 (a)\le n\omega_n\int_0^aG(r)\frac{dr}{r}$, and by L'H\^opital's rule
\[
\lim_{a\to\infty}\frac{\int_0^aG(r)\frac{dr}{r}}{G(a)}=\lim_{a\to\infty}\frac{G(a)}{a\,g(a)}=\lim_{a\to\infty}\frac{e^{a^2}-1}{2a^2e^{a^2}}=0,
\]
so that $ G_0 (a)/G(a)\to0$ as $a\to\infty$. The equivalence \eqref{equiv.escala}, with scaling of the argument, is the natural one in Orlicz spaces.
\end{rem}

By \eqref{exacto} (or directly from the definition, by Fubini's theorem), $ G_0 $ can be written as
\begin{equation}\label{Gt.F}
 G_0 (a)=\int_{\mathbb S^{n-1}}F(a|\omega\cdot e_n|)\,d\mathcal H^{n-1},\qquad F(b)=\int_0^bG(r)\,\frac{dr}{r}.
\end{equation}
Moreover, if $f$ is integrable on $[-1,1]$, then 
\[\int_{\mathbb S^{n-1}}f(\omega\cdot e_n)\,d\mathcal H^{n-1}=\mathcal H^{n-2}(\mathbb S^{n-2})\int_{-1}^1f(\tau)(1-\tau^2)^{\frac{n-3}2}\,d\tau\] 
for $n\ge2$. In particular, for $p>-1$,
\[
K_{n,p}:=\int_{\mathbb S^{n-1}}|\omega\cdot e_n|^p\,d\mathcal H^{n-1}=\frac{2\pi^{\frac{n-1}2}\,\Gamma(\frac{p+1}2)}{\Gamma(\frac{n+p}2)},
\]
and for $n=3$ the weight is constant: $\int_{\mathbb S^2}f(\omega\cdot e_3)\,d\mathcal H^2=2\pi\int_{-1}^1f(\tau)\,d\tau$.

\begin{exam}
The following examples correct and simplify those in \cite[Section 2.4]{FBS}.
\begin{enumerate}
\item[(i)] If $G(a)=a^p$ with $p>1$, then $F(b)=\frac{b^p}p$ and $ G_0 (a)=\frac{K_{n,p}}{p}\,a^p$.
\item[(ii)] If $G=\sum_kc_kG_k$ with $c_k\ge0$, then $ G_0 =\sum_kc_k(G_k)_0$.
\item[(iii)] If $G(a)=\max\{a^p,a^q\}$ with $1<q<p$, then $F(b)=\frac{b^q}q$ for $b\le1$ and $F(b)=\frac{b^p}p+\frac1q-\frac1p$ for $b>1$. Hence $ G_0 (a)=\frac{K_{n,q}}q\,a^q$ for $a\le1$, while for $a>1$
\begin{align*}
 G_0 (a)={}&\frac{a^q}q\int_{\{|\omega\cdot e_n|\le\frac1a\}}|\omega\cdot e_n|^q\,d\mathcal H^{n-1}+\frac{a^p}p\int_{\{|\omega\cdot e_n|>\frac1a\}}|\omega\cdot e_n|^p\,d\mathcal H^{n-1}\\
&+\Big(\frac1q-\frac1p\Big)\mathcal H^{n-1}\Big(\Big\{|\omega\cdot e_n|>\frac1a\Big\}\Big).
\end{align*}
For $n=3$ this gives, for $a>1$,
\[
 G_0 (a)=4\pi\Big(\frac{a^{-1}}{q(q+1)}+\frac{a^p-a^{-1}}{p(p+1)}+\Big(\frac1q-\frac1p\Big)\Big(1-\frac1a\Big)\Big).
\]
\item[(iv)] If $G(a)=a^p\log^\alpha(1+a)$ with $p>1$ and $\alpha>0$, then $G$ is a Young function satisfying the $\Delta_2$ condition, and $ G_0 $ is not elementary. However, from \eqref{Gt.F} one easily obtains its behavior at $0$ and at $\infty$:
\[
 G_0 (a)\sim\frac{K_{n,p+\alpha}}{p+\alpha}\,a^{p+\alpha}\quad\text{as }a\to0,\qquad G_0 (a)\sim\frac{K_{n,p}}{p}\,a^p\log^\alpha a\quad\text{as }a\to\infty.
\]
\end{enumerate}
Typical examples of Young functions not satisfying $\Delta_2$, to which the results of this paper apply, are $G(a)=e^{a^2}-1$, $G(a)=e^a-1-a$ and $G(a)=e^{a^\gamma}-1$ for $\gamma>1$.
\end{exam}

\begin{exam}
For some Young functions without the $\Delta_2$ condition, $ G_0 $ can be computed explicitly from \eqref{Gt.F}. For $n=1$ we have $ G_0 (a)=2F(a)$, and for $n=3$, an integration by parts gives $ G_0 (a)=4\pi\big(F(a)-\frac1a\int_0^aG(r)\,dr\big)$. Let $\operatorname{Ei}(x)=-\operatorname{p.v.}\int_{-x}^\infty\frac{e^{-r}}{r}\,dr$ for $x>0$ be the exponential integral, so that $\int_0^x\frac{e^r-1}{r}\,dr=\operatorname{Ei}(x)-\gamma_{\rm E}-\log x$, where $\gamma_{\rm E}$ is the Euler--Mascheroni constant, and let $\operatorname{erfi}(x)=\frac2{\sqrt\pi}\int_0^xe^{r^2}\,dr$.
\begin{itemize}
\item[(i)] For $G(a)=e^a-1-a$,
\begin{align*}
& G_0 (a)=2\big(\operatorname{Ei}(a)-\gamma_{\rm E}-\log a-a\big)\quad(n=1),\\
& G_0 (a)=4\pi\Big(\operatorname{Ei}(a)-\gamma_{\rm E}-\log a-\frac{e^a-1-a}{a}\Big)-2\pi a\quad(n=3).
\end{align*}
\item[(ii)] For $G(a)=e^{a^2}-1$,
\begin{align*}
& G_0 (a)=\operatorname{Ei}(a^2)-\gamma_{\rm E}-2\log a\quad(n=1),\\
& G_0 (a)=2\pi\big(\operatorname{Ei}(a^2)-\gamma_{\rm E}-2\log a\big)-\frac{2\pi^{3/2}\operatorname{erfi}(a)}{a}+4\pi\quad(n=3).
\end{align*}
\end{itemize}
All formulas in this example that contain $\log a$ or a quotient by $a$ are written for $a>0$; their values at $a=0$ are understood by continuous extension, with $ G_0 (0)=0$.
Since $\operatorname{Ei}(x)\sim\frac{e^x}{x}$ as $x\to\infty$, for $n=1$ we get $ G_0 (a)\sim\frac{2e^a}{a}$ in the first case and $ G_0 (a)\sim\frac{e^{a^2}}{a^2}$ in the second one. In both cases $ G_0 (a)/G(a)\to0$ as $a\to\infty$, in accordance with the previous remark.

More generally, if $G(a)=\sum_{k\ge2}c_ka^k$ with $c_k\ge0$, integrating term by term, we get that  $ G_0 (a)=\sum_{k\ge2}\frac{c_kK_{n,k}}{k}\,a^k$.
\end{exam}

\subsection{A pointwise estimate}

For $v$ measurable and $0<s<1$ we consider the pointwise density of the modular
\[
\mathcal G_s[v](x)=(1-s)\int_{\mathbb{R}^n} G\Big(\frac{|v(x)-v(y)|}{|x-y|^s}\Big)\frac{dy}{|x-y|^n},\qquad\text{so that}\qquad \Phi_{G,\nu_s}(D^sv)=\int_{\mathbb{R}^n}\mathcal G_s[v]\,dx .
\]
For $v\in C^2(\mathbb{R}^n)$ we write $\|v\|_{C^2}=\|v\|_\infty+\|\nabla v\|_\infty+\|D^2v\|_\infty$.

\begin{lema}\label{lema.puntual}
Let $K>0$ and $v\in C^2(\mathbb{R}^n)$ with $\|v\|_{C^2}\le K$. Then, for every $x\in\mathbb{R}^n$ and $s\in(0,1)$,
\[
\big|\mathcal G_s[v](x)- G_0 (|\nabla v(x)|)\big|\le(1-s)\,n\omega_n\Big(\frac{L_K K}{2}+\frac{G(2K)}{s}\Big),
\]
where $L_K$ is the Lipschitz constant of $G$ on $[0,K]$. In particular, $\mathcal G_s[v]\to G_0 (|\nabla v|)$ uniformly in $\mathbb{R}^n$ and uniformly for $v$ in bounded subsets of $C^2(\mathbb{R}^n)$.
\end{lema}

\begin{proof}
For $|h|<1$, both $\frac{|v(x+h)-v(x)|}{|h|^s}$ and $\frac{|\nabla v(x)\cdot h|}{|h|^s}$ belong to $[0,K]$, and their difference is at most $\frac K2|h|^{2-s}$. Hence
\[
(1-s)\int_{|h|<1}\Big|G\Big(\frac{|v(x+h)-v(x)|}{|h|^s}\Big)-G\Big(\frac{|\nabla v(x)\cdot h|}{|h|^s}\Big)\Big|\frac{dh}{|h|^n}\le(1-s)\frac{L_KK}{2}\,\frac{n\omega_n}{2-s}.
\]
By polar coordinates, rotation invariance and \eqref{exacto},
\[
\begin{aligned}
(1-s)\int_{|h|<1}G\Big(\frac{|\nabla v(x)\cdot h|}{|h|^s}\Big)\frac{dh}{|h|^n}&=(1-s)\int_0^1\int_{\mathbb S^{n-1}}G\big(|\nabla v(x)||\omega\cdot e_n|r^{1-s}\big)\,d\mathcal H^{n-1}\,\frac{dr}{r}\\&= G_0 (|\nabla v(x)|).
\end{aligned}
\]
Finally, for $|h|\ge1$, \eqref{P5} gives $G(|v(x+h)-v(x)|/|h|^s)\le|h|^{-s}G(2K)$, and the contribution of $\{|h|\ge1\}$ is at most $(1-s)\frac{n\omega_n}{s}G(2K)$.
\end{proof}

\subsection{An upper bound}

\begin{lema}\label{lema.cota}
Let $u\in L_G(\mathbb{R}^n)$ be weakly differentiable with $\Phi_{ G_0 }(|\nabla u|)<\infty$. Then, for every $s\in(0,1)$,
\[
\Phi_{G,\nu_s}(D^su)\le\Phi_{ G_0 }(|\nabla u|)+\frac{1-s}{s}\,n\omega_n\,\Phi_G(2u).
\]
In particular, $W^1L_G(\mathbb{R}^n)\subset W^sL_G(\mathbb{R}^n)$ and $W^1E_G(\mathbb{R}^n)\subset W^sE_G(\mathbb{R}^n)$ for all $s\in(0,1)$, and
\[
[u]_{s,G}\le\Big(1+n\omega_n\frac{1-s}{s}\Big)\big([u]_{1, G_0 }+2\|u\|_G\big)\le3\beta_n\Big(1+n\omega_n\frac{1-s}{s}\Big)\|u\|_{1,G}
\]
for every $u\in W^1L_G(\mathbb{R}^n)$.
\end{lema}

\begin{proof}
This is \cite[Lemma 3.6]{CFBS}. For $|h|<1$ write $h=r\omega$ and $u(x+h)-u(x)=\int_0^1\nabla u(x+\theta h)\cdot h\,d\theta$ for a.e.\ $x$. Jensen's inequality and translation invariance give
\[
\int_{\mathbb{R}^n} G\Big(\frac{|u(x+h)-u(x)|}{|h|^s}\Big)dx\le\int_{\mathbb{R}^n} G\big(r^{1-s}|\nabla u(x)\cdot\omega|\big)\,dx ,
\]
and integrating against $(1-s)\,d\mathcal H^{n-1}(\omega)\frac{dr}{r}$, by \eqref{exacto}, the contribution of $\{|h|<1\}$ is at most $\Phi_{ G_0 }(|\nabla u|)$. For $|h|\ge1$ use \eqref{P5} and \eqref{P4}: $G(|u(x+h)-u(x)|/|h|^s)\le\frac{|h|^{-s}}2(G(2|u(x+h)|)+G(2|u(x)|))$.

For the norm inequality, apply the estimate to $u/\lambda$ with $\lambda=[u]_{1, G_0 }+2\|u\|_G$. Since $\Phi_{ G_0 }(|\nabla u|/\lambda)\le1$ and $\Phi_G(2u/\lambda)\le1$, we get
\[
\Phi_{G,\nu_s}\Big(\frac{D^su}{\lambda}\Big)\le1+n\omega_n\frac{1-s}{s},
\]
so $[u]_{s,G}\le(1+n\omega_n\frac{1-s}{s})\lambda$ by \eqref{mod-norm}, and the last inequality follows from \eqref{equiv.normas.1}, since $\beta_n\ge1$. The inclusion $W^1E_G\subset W^sE_G$ follows in the same way, applying the estimate to $ku$ for every $k>0$ and using \eqref{equiv.escala}.
\end{proof}

\section{Compactness}\label{sec.comp}

\begin{thm}\label{teo.comp}
Let $0<s<1$ and let $\mathcal K\subset W^sL_G(\mathbb{R}^n)$ be bounded. Then for every $R>0$, the set $\{u|_{B_R}\colon u\in\mathcal K\}$ is contained in $E_G(B_R)$ and is precompact in $L_G(B_R)$. In particular, $W^sL_G(\mathbb{R}^n)\subset E_{G,{\rm loc}}(\mathbb{R}^n)$, and if $\Omega$ is bounded then the inclusion $W^s_0L_G(\Omega)\subset E_G(\Omega)$ is compact.
\end{thm}

\begin{proof}
Let $M=\sup_{u\in\mathcal K}\|u\|_{s,G}$. For $u\in\mathcal K$ and $\varepsilon\in(0,1]$, let $\lambda=\sup_{|h|\le\varepsilon}\|\tau_hu-u\|_G\le C_nM\varepsilon^s$, by Lemma \ref{lema.translation}. Since $u_\varepsilon-u=\int(\tau_{-y}u-u)\rho_\varepsilon(y)\,dy$, Jensen's inequality gives
\[
\int_{\mathbb{R}^n} G\Big(\frac{|u_\varepsilon-u|}{\lambda}\Big)dx\le\int\Big(\int_{\mathbb{R}^n} G\Big(\frac{|\tau_{-y}u-u|}{\lambda}\Big)dx\Big)\rho_\varepsilon(y)\,dy\le1,
\]
so $\|u_\varepsilon-u\|_G\le C_nM\varepsilon^s$ uniformly in $\mathcal K$. 

By Lemma \ref{lema.C2}, $\{u_\varepsilon|_{\overline B_R}\colon u\in\mathcal K\}$ is bounded in $C^1(\overline B_R)$, hence precompact in $C(\overline B_R)$ by the Arzel\`a--Ascoli theorem, and therefore precompact in $L_G(B_R)$, since $\|v\|_{L_G(B_R)}\le C_R\|v\|_\infty$. Thus $\mathcal K|_{B_R}$ is totally bounded in $L_G(B_R)$.

Moreover, every $u|_{B_R}$ is at distance at most $C_nM\varepsilon^s$ from the bounded function $u_\varepsilon|_{B_R}$, for every $\varepsilon$; since $E_G(B_R)$ is the closure of the bounded functions, $u|_{B_R}\in E_G(B_R)$. If $\Omega\subset B_R$, functions in $W^s_0L_G(\Omega)$ vanish outside $B_R$, which gives the last assertion.
\end{proof}

\begin{rem}
Theorem \ref{teo.comp} corrects and extends \cite[Theorem 3.1]{FBS}, whose proof used the $\Delta_2$ condition through a modular version of Lemma \ref{lema.translation}. Without $\Delta_2$, a modular estimate of $\tau_hu-u$ does not give a norm estimate, and the norm version of Lemma \ref{lema.translation} is needed. See also \cite[Theorem B.3]{FBSp}.
\end{rem}

\section{The Bourgain--Brezis--Mironescu formula}\label{sec.bbm}

The following lemma is the closure argument used in \cite[proof of Theorem 3.5]{CFBS}.

\begin{lema}\label{lema.cierre}
Let $u\in L_G(\mathbb{R}^n)$, and let $\{v_j\}\subset W^{1,1}_{\rm loc}(\mathbb{R}^n)$ be such that $v_j\to u$ in $L^1_{\rm loc}(\mathbb{R}^n)$ and $\Phi_{ G_0 }(|\nabla v_j|)\le M$ for all $j$. Then $u\in W^1L_G(\mathbb{R}^n)$ and
\[
\Phi_{ G_0 }(|\nabla u|)\le\liminf_{j\to\infty}\Phi_{ G_0 }(|\nabla v_j|).
\]
\end{lema}

\begin{proof}
By Proposition \ref{prop.equiv}, $\int_{\mathbb{R}^n} G(|\nabla v_j|/4)\le M/c_n$. Since $G$ is superlinear at infinity, the de la Vall\'ee-Poussin criterion shows that $\{\nabla v_j\}$ is equi-integrable on bounded sets, and by the Dunford--Pettis theorem a subsequence converges weakly in $L^1_{\rm loc}$ to some $V$. Since $v_j\to u$ in $L^1_{\rm loc}$, we get that  $V=\nabla u$. 

The functionals $w\mapsto\int G(|w|/4)$ and $w\mapsto\int G_0 (|w|)$ are convex and strongly lower semicontinuous, hence weakly lower semicontinuous in $L^1_{\rm loc}$. Therefore $\int G(|\nabla u|/4)\le M/c_n$, so $|\nabla u|\in L_G$, and the inequality holds (along a subsequence realizing the $\liminf$).
\end{proof}

\begin{thm}\label{main1}
Let $G$ be a Young function and $u\in L_G(\mathbb{R}^n)$. Then
\begin{equation}\label{bbm.liminf}
\Phi_{ G_0 }(|\nabla u|)\le\liminf_{s\uparrow1}\Phi_{G,\nu_s}(D^su),
\end{equation}
where the left-hand side is $+\infty$ if $u\notin W^1L_G(\mathbb{R}^n)$. If moreover $\Phi_G(2u)<\infty$ (for instance, if $u\in E_G(\mathbb{R}^n)$), then
\begin{equation}\label{bbm}
\lim_{s\uparrow1}\Phi_{G,\nu_s}(D^su)=\Phi_{ G_0 }(|\nabla u|).
\end{equation}
\end{thm}

\begin{proof}
Let $L$ be the right-hand side of \eqref{bbm.liminf}; we may assume $L<\infty$ and take $s_k\uparrow1$ with $\Phi_{G,\nu_{s_k}}(D^{s_k}u)\to L$. Fix $\varepsilon>0$. By Lemma \ref{lema.C2}, $u_\varepsilon\in C^2$ with bounded derivatives up to order 2. By Lemma \ref{lema.reg}, Fatou's lemma and Lemma \ref{lema.puntual},
\[
L\ge\liminf_{k}\Phi_{G,\nu_{s_k}}(D^{s_k}u_\varepsilon)\ge\int_{\mathbb{R}^n}\liminf_k\mathcal G_{s_k}[u_\varepsilon]\,dx=\Phi_{ G_0 }(|\nabla u_\varepsilon|).
\]
Since $u_\varepsilon\to u$ in $L^1_{\rm loc}$ as $\varepsilon\to0$, Lemma \ref{lema.cierre} gives $u\in W^1L_G(\mathbb{R}^n)$ and $\Phi_{ G_0 }(|\nabla u|)\le L$.

If $\Phi_G(2u)<\infty$, we may assume $\Phi_{ G_0 }(|\nabla u|)<\infty$, and Lemma \ref{lema.cota} gives
\[
\limsup_{s\uparrow1}\Phi_{G,\nu_s}(D^su)\le\Phi_{ G_0 }(|\nabla u|).
\]
\end{proof}

\begin{rem}
The proof of \cite[Theorem 4.1]{FBS} extends the result from $C^2_c$ to $W^{1,G}$ by density in norm, which requires $\Delta_2$. The argument above avoids density. The equality \eqref{bbm} holds in particular for every $u\in W^1E_G(\mathbb{R}^n)$, and the condition $\Phi_G(2u)<\infty$ is only used to control the contribution of $\{|h|\ge1\}$. Compare with \cite[Theorem 1.1]{ACPS}.
\end{rem}

\section{The case of a sequence}\label{sec.seq}

\subsection{A uniform estimate}
The following result replaces \cite[Theorem 5.2]{FBS}, whose proof relies on \cite[Lemma 5.3]{FBS}, which does not hold as stated (see Remark \ref{rem.lema53}).

\begin{thm}\label{teo.s1.s2}
Let $G$ be a Young function and $u\in L_G(\mathbb{R}^n)$. Then, for every $0<t\le s<1$,
\begin{equation}\label{s1s2.mod}
\Phi_{G,\nu_t}\Big(\frac{D^tu}8\Big)\le3\,\Phi_{G,\nu_s}(D^su)+n\omega_n\,\frac{1-t}{t}\,\Phi_G\Big(\frac u4\Big),
\end{equation}
where the modular inequality is understood with extended nonnegative values. In particular, if $u\in W^sL_G(\mathbb{R}^n)$, then
\begin{equation}\label{s1s2.norm}
[u]_{t,G}\le C_t\|u\|_{s,G},\qquad\text{and hence}\qquad\|u\|_{t,G}\le(1+C_t)\|u\|_{s,G},\qquad C_t=8\Big(3+n\omega_n\frac{1-t}{t}\Big).
\end{equation}
\end{thm}

\begin{proof}
For $r>0$ and $\varrho\ge0$ let
\[
\Psi(r,\varrho)=\int_{\mathbb S^{n-1}}\int_{\mathbb{R}^n} G\Big(\varrho\,\frac{|u(x+r\omega)-u(x)|}{r}\Big)dx\,d\mathcal H^{n-1}(\omega).
\]
$\Psi(r,\cdot)$ is nondecreasing and, since $\frac{|u(x+2r\omega)-u(x)|}{2r}$ is bounded by the average of the corresponding quotients at scale $r$ at the points $x$ and $x+r\omega$, convexity and translation invariance give
\begin{equation}\label{doubling}
\Psi(2r,\varrho)\le\Psi(r,\varrho)\qquad\text{for all }r>0,\ \varrho\ge0 .
\end{equation}
Let $I_k=[2^{-k-1},2^{-k}]$ and $\varphi_k(\varrho)=\int_{I_k}\Psi(r,\varrho)\frac{dr}{r}$, $k\ge0$. Changing variables $r=2^{-(j-k)}r'$ 
\begin{equation}\label{monot}
\varphi_j(\varrho)\ge\varphi_k(\varrho)\qquad\text{for }j\ge k .
\end{equation}
where we have used \eqref{doubling}.

Let $\alpha=1-t$, $\beta=1-s$ (so $0<\beta\le\alpha\le1$) and $\lambda_0=\frac18$. In polar coordinates,
\[
N_t:=\int_{|h|<1}\int_{\mathbb{R}^n} G\Big(\frac{|u(x+h)-u(x)|}{8|h|^t}\Big)dx\frac{dh}{|h|^n}=\int_0^1\Psi(r,\lambda_0 r^\alpha)\frac{dr}{r}\le\sum_{k\ge0}a_k,
\]
where $a_k=\varphi_k(\lambda_02^{-k\alpha})$, and
\[
N_s:=\int_{|h|<1}\int_{\mathbb{R}^n} G\Big(\frac{|u(x+h)-u(x)|}{|h|^s}\Big)dx\frac{dh}{|h|^n}=\int_0^1\Psi(r,r^\beta)\frac{dr}{r}\ge\sum_{j\ge0}b_j,
\]
where $b_j=\varphi_j(2^{-(j+1)\beta})$.
Let $M=\lceil1/\beta\rceil$, so that $1\le M\beta\le2$, and for $k\ge0$ let \[
B_k=\{j\colon\sigma_k\le j\le\sigma_k+M-1\}\]
with $\sigma_k=\lceil k\alpha/\beta\rceil$. If $j\in B_k$, then $j\ge k$ and $(j+1)\beta\le k\alpha+3$, so $2^{-(j+1)\beta}\ge\lambda_02^{-k\alpha}$; by the monotonicity of $\varphi_j$ and \eqref{monot}, $b_j\ge\varphi_j(\lambda_02^{-k\alpha})\ge a_k$. Each $B_k$ has $M\ge1/\beta$ elements, and each $j$ belongs to at most $\frac2\alpha+1\le\frac3\alpha$ of the sets $B_k$. Hence
\[
\frac1\beta\sum_ka_k\le\sum_k\sum_{j\in B_k}b_j\le\frac3\alpha\sum_jb_j,\qquad\text{that is,}\qquad N_t\le\frac{3(1-s)}{1-t}N_s\le\frac{3(1-s)}{1-t}\Phi_{G,\nu}(D^su).
\]
For $|h|\ge1$, \eqref{P5} and \eqref{P4} give \[
G\big(\frac{|u(x+h)-u(x)|}{8|h|^t}\big)\le\frac{|h|^{-t}}2\big(G(\frac{|u(x+h)|}4)+G(\frac{|u(x)|}4)\big),
\] 
and the contribution of $\{|h|\ge1\}$ to $\Phi_{G,\nu}(D^tu/8)$ is at most $\frac{n\omega_n}{t}\Phi_G(u/4)$. Multiplying by $1-t$ we obtain \eqref{s1s2.mod}.

To prove \eqref{s1s2.norm}, let $\lambda=\|u\|_{s,G}$. Applying \eqref{s1s2.mod} to $u/\lambda$, and since $\Phi_{G,\nu_s}(D^su/\lambda)\le1$ and $\Phi_G(u/(4\lambda))\le1$, we get $\Phi_{G,\nu_t}\big(\frac{D^tu}{8\lambda}\big)\le3+n\omega_n\frac{1-t}{t}$, and \eqref{mod-norm} gives the result.
\end{proof}

\begin{rem}\label{rem.lema53}
\cite[Lemma 5.3]{FBS} states that if $g,h\colon(0,1)\to\mathbb{R}^+$, $g(\sigma)\le c\,g(\sigma/2)$ and $h$ is decreasing, then \[
\int_0^1\sigma^rgh\,d\sigma\ge\frac{r+1}{2^{r+1}}\int_0^1\sigma^rg\,d\sigma\int_0^1\sigma^rh\,d\sigma \text{ for } r>-1.
\]
This is true for $c=1$ (see \cite{BBM}), but false for $c>1$, even with any other positive constant: for $g(\sigma)=\sigma^m$ (so $c=2^m$) and $h=\chi_{(0,\delta)}$, the ratio between both sides is a constant times $\delta^m\to0$. The proof of Theorem \ref{teo.s1.s2} above uses only the case $c=1$, in the form \eqref{doubling}.
\end{rem}

\subsection{Compactness and lower semicontinuity}

\begin{thm}\label{main2}
Let $G$ be a Young function, $s_k\uparrow1$, and $\{u_k\}\subset L_G(\mathbb{R}^n)$ with
\[
\sup_k\Phi_{G,\nu_{s_k}}(D^{s_k}u_k)<\infty,\qquad\sup_k\Phi_G(u_k)<\infty .
\]
Then there exist $u\in W^1L_G(\mathbb{R}^n)\cap E_{G,{\rm loc}}(\mathbb{R}^n)$ and a subsequence such that $u_{k_j}\to u$ in $L_{G,{\rm loc}}(\mathbb{R}^n)$, and
\[
\Phi_{ G_0 }(|\nabla u|)\le\liminf_{k\to\infty}\Phi_{G,\nu_{s_k}}(D^{s_k}u_k).
\]
\end{thm}

\begin{proof}
Fix $t\in(0,1)$. By \eqref{mod-norm}, $\sup_k\|u_k\|_{s_k,G}<\infty$, so by \eqref{s1s2.norm}, $\{u_k\}_{s_k\ge t}$ is bounded in $W^tL_G(\mathbb{R}^n)$. By Theorem \ref{teo.comp} and a diagonal argument, a subsequence converges in $L_G(B_R)$ for every $R$ to some $u\in E_{G,{\rm loc}}(\mathbb{R}^n)$; by Fatou's lemma, $u\in L_G(\mathbb{R}^n)$. The lower bound follows from Theorem \ref{teo.liminf} below.
\end{proof}

\begin{thm}\label{teo.liminf}
Let $s_k\uparrow1$ and let $\{u_k\}\subset L_G(\mathbb{R}^n)$ be bounded in $L_G(\mathbb{R}^n)$ and such that $u_k\to u$ in $L^1_{\rm loc}(\mathbb{R}^n)$. Then $u\in L_G(\mathbb{R}^n)$ and
\[
\Phi_{ G_0 }(|\nabla u|)\le\liminf_{k\to\infty}\Phi_{G,\nu_{s_k}}(D^{s_k}u_k),
\]
where the left-hand side is $+\infty$ if $u\notin W^1L_G(\mathbb{R}^n)$.
\end{thm}

\begin{proof}
Let $L$ be the right-hand side; we may assume $L<\infty$ and, passing to a subsequence, that $\Phi_{G,\nu_{s_k}}(D^{s_k}u_k)\to L$. Since $u_k\to u$ in $L^1_{\rm loc}$, after passing to a further subsequence we have $u_k(x)\to u(x)$ for almost every $x\in\mathbb{R}^n$. Boundedness in $L_G$ gives a constant $C>0$ such that $\Phi_G(u_k/C)\le1$ for every $k$; Fatou's lemma then yields $\Phi_G(u/C)\le1$, and therefore $u\in L_G(\mathbb{R}^n)$. Fix $\varepsilon>0$ and $R>0$. By Lemma \ref{lema.C2}, $\{(u_k)_\varepsilon\}_k$ is bounded in $C^2(\mathbb{R}^n)$ and $(u_k)_\varepsilon\to u_\varepsilon$ in $C^1(\overline B_R)$. By Lemma \ref{lema.reg} and Lemma \ref{lema.puntual},
\[
\Phi_{G,\nu_{s_k}}(D^{s_k}u_k)\ge\int_{B_R}\mathcal G_{s_k}[(u_k)_\varepsilon]\,dx\ge\int_{B_R} G_0 (|\nabla(u_k)_\varepsilon|)\,dx-|B_R|\,\eta(1-s_k),
\]
with $\eta(\sigma)\to0$ as $\sigma\to0$. Letting $k\to\infty$ and using the continuity of $ G_0 $, $\int_{B_R} G_0 (|\nabla u_\varepsilon|)\le L$. Letting $R\to\infty$, $\Phi_{ G_0 }(|\nabla u_\varepsilon|)\le L$, and Lemma \ref{lema.cierre} concludes the proof.
\end{proof}

\begin{rem}
This proof does not use truncations, so it also simplifies the proof of \cite[Theorem 3.5]{CFBS}. Note that only convergence in $L^1_{\rm loc}$ and boundedness in $L_G$ are needed.
\end{rem}

\section{Consequences and applications}\label{sec.app}

\subsection{Poincar\'e inequalities}
We state in this section some important consequences which follow from our main results.
\begin{prop}\label{prop.poinc.s}
Let $\Omega\subset\mathbb{R}^n$ be open and bounded, $d=\operatorname{diam}\Omega$, and let $u$ be measurable with $u=0$ a.e.\ in $\mathbb{R}^n\setminus\Omega$. Then, for every $s\in(0,1)$,
\begin{equation}\label{poinc.mod.s}
n\omega_n\log2\int_\Omega G\Big(\frac{|u|}{16(d+1)}\Big)dx\le\Phi_{G,\nu_s}(D^su),
\end{equation}
and, if $u\in W^1_0L_G(\Omega)$,
\begin{equation}\label{poinc.mod.1}
\int_\Omega G\Big(\frac{|u|}{d}\Big)dx\le\int_\Omega G(|\nabla u|)\,dx,\qquad\int_\Omega G\Big(\frac{\alpha_n|u|}{d}\Big)dx\le\Phi_{ G_0 }(|\nabla u|).
\end{equation}
In particular, letting $R_\Omega=\max\{16(d+1),d/\alpha_n\}$ and $\tilde c_n=\min\{1,n\omega_n\log2\}$,
\begin{equation}\label{poinc.mod}
\tilde c_n\int_\Omega G\Big(\frac{|u|}{R_\Omega}\Big)dx\le\mathcal J_s(u)\qquad\text{for every }s\in(0,1]\text{ and }u\in W^s_0L_G(\Omega),
\end{equation}
where $\mathcal J_s(u)=\Phi_{G,\nu_s}(D^su)$ for $s<1$ and $\mathcal J_1(u)=\Phi_{ G_0 }(|\nabla u|)$.
\end{prop}

\begin{proof}
\emph{The case $s<1$.} We use the notation of Step 1 in the proof of Lemma \ref{lema.translation}, with $r_0=2(d+1)$. The argument there, which only uses the monotonicity $\Psi(2r,\varrho)\le\Psi(r,\varrho)$ and not the normalization, gives
\[
M\varphi_0\Big(\frac{r_0^{1-s}}8\Big)\le\int_0^{r_0}\Psi(r,r^{1-s})\frac{dr}r\le\Phi_{G,\nu}(D^su),\qquad M=\Big\lceil\frac1{1-s}\Big\rceil,
\]
and since $M(1-s)\ge1$, we get $\varphi_0(r_0^{1-s}/8)\le\Phi_{G,\nu_s}(D^su)$. If $x\in\Omega$ and $r\ge d+1$, then $x+r\omega\notin\Omega$, so $u(x+r\omega)=0$. Hence, for $\varrho>0$,
\begin{align*}
\varphi_0(\varrho)=\int_{d+1}^{2(d+1)}\Psi(r,\varrho)\frac{dr}r&\ge\int_{d+1}^{2(d+1)}\int_{\mathbb S^{n-1}}\int_\Omega G\Big(\frac{\varrho|u(x)|}{r}\Big)dx\,d\mathcal H^{n-1}\,\frac{dr}r\\
&\ge n\omega_n\log2\int_\Omega G\Big(\frac{\varrho|u|}{r_0}\Big)dx .
\end{align*}
With $\varrho=r_0^{1-s}/8$ we have $\varrho/r_0=r_0^{-s}/8\ge1/(8r_0)$, since $r_0\ge1$, and \eqref{poinc.mod.s} follows.

\emph{The case $s=1$.} Translating, we may assume $\Omega\subset\{0<x_1<d\}$. Since $u\in W^{1,1}_{\rm loc}(\mathbb{R}^n)$ vanishes outside $\Omega$, for a.e.\ $x'\in\mathbb{R}^{n-1}$ the function $y_1\mapsto u(y_1,x')$ is absolutely continuous and vanishes for $y_1\le0$, so $|u(x)|\le\int_0^d|\nabla u(y_1,x')|\,dy_1$ for a.e.\ $x\in\Omega$. By Jensen's inequality,
\[
G\Big(\frac{|u(x)|}{d}\Big)\le\frac1d\int_0^dG(|\nabla u(y_1,x')|)\,dy_1,
\]
and integrating over $x\in\{0<x_1<d\}$ we get the first inequality in \eqref{poinc.mod.1}. Here we use that $\nabla u=0$ almost everywhere on $\mathbb{R}^n\setminus\Omega$: a Sobolev function that is zero almost everywhere on a measurable set has zero weak gradient almost everywhere on that set. Thus the integral of $G(|\nabla u|)$ over the strip reduces to the integral over $\Omega$. Applying it to $\alpha_nu$ and using \eqref{equiv.escala} we get the second one.

Finally, \eqref{poinc.mod} follows from \eqref{poinc.mod.s}, \eqref{poinc.mod.1} and the monotonicity of $G$.
\end{proof}

Proposition \ref{prop.poinc.s} replaces \cite[Theorem 6.1 and Corollary 6.2]{FBS}. There, a modular Poincar\'e inequality uniform in $s$ was proved under the $\Delta_2$ condition, by contradiction and compactness, following \cite{Ponce}. Without $\Delta_2$, the inequality must include a scaling of the argument, as in \eqref{poinc.mod}, and the proof above is direct and gives explicit constants. As a consequence we obtain the following norm version.

\begin{thm}\label{teo.poinc.unif}
Let $\Omega\subset\mathbb{R}^n$ be open and bounded. Then
\[
\|u\|_{L_G(\Omega)}\le C_\Omega\,[u]_{s,G}\qquad\text{for every }s\in(0,1]\text{ and every }u\in W^s_0L_G(\Omega),
\]
where $C_\Omega=R_\Omega/\tilde c_n$, with $R_\Omega$ and $\tilde c_n$ as in Proposition \ref{prop.poinc.s}.
\end{thm}

\begin{proof}
For $s=1$, the first inequality in \eqref{poinc.mod.1} applied to $u/[u]_{1,G}$ gives \[\|u\|_{L_G(\Omega)}\le d\,[u]_{1,G}\le C_\Omega[u]_{1,G}.\] 
Let $s<1$. We may assume $\lambda=[u]_{s,G}>0$. Then $\mathcal J_s(u/\lambda)\le1$, and \eqref{poinc.mod} gives $\int_\Omega G\big(\frac{|u|}{R_\Omega\lambda}\big)dx\le\frac1{\tilde c_n}$. Since $\tilde c_n\le1$, the convexity of $G$ gives $\int_\Omega G\big(\frac{\tilde c_n|u|}{R_\Omega\lambda}\big)dx\le1$, that is, $\|u\|_{L_G(\Omega)}\le\frac{R_\Omega}{\tilde c_n}\lambda$.
\end{proof}

\subsection{$\Gamma$-convergence}

Let $\Omega\subset\mathbb{R}^n$ be open. On $E_G(\Omega)$, endowed with the norm topology, define
\[
\mathcal J_s(u)=\begin{cases}\Phi_{G,\nu_s}(D^su)&u\in W^s_0L_G(\Omega),\\+\infty&\text{otherwise},\end{cases}\qquad\mathcal J(u)=\begin{cases}\Phi_{ G_0 }(|\nabla u|)&u\in W^1_0L_G(\Omega),\\+\infty&\text{otherwise}.\end{cases}
\]

\begin{thm}\label{teo.gamma.conv}
For every open set $\Omega\subset\mathbb{R}^n$, $\mathcal J=\Gamma\text{-}\lim_{s\uparrow1}\mathcal J_s$ in $E_G(\Omega)$. Moreover, the constant sequence is a recovery sequence.
\end{thm}

\begin{proof}
\emph{Liminf.} Let $u_k\to u$ in $E_G(\Omega)$ and $s_k\uparrow1$, with $\liminf\mathcal J_{s_k}(u_k)<\infty$. Along a subsequence, $u_k\in W^{s_k}_0L_G(\Omega)$ vanish outside $\Omega$, are bounded in $L_G$ and converge in $L^1_{\rm loc}$. By Theorem \ref{teo.liminf}, $u\in W^1L_G(\mathbb{R}^n)$, $u=0$ a.e.\ outside $\Omega$, and $\mathcal J(u)\le\liminf\mathcal J_{s_k}(u_k)$.

\emph{Limsup.} Let $u\in E_G(\Omega)$ with $\mathcal J(u)<\infty$. Since $\Phi_G(2u)<\infty$, Lemma \ref{lema.cota} shows that $u\in W^s_0L_G(\Omega)$ for every $s$ and that $\limsup_{s\uparrow1}\mathcal J_s(u)\le\mathcal J(u)$.
\end{proof}

\begin{rem}\label{rem.ae}
\begin{enumerate}
\item[(i)] The effective domain of $\mathcal J$ is the modular domain \[\{u\in W^1_0L_G(\Omega)\cap E_G(\Omega)\colon\Phi_{ G_0 }(|\nabla u|)<\infty\},\] which in general is strictly larger than $W^1_0E_G(\Omega)$ (for instance, if $\Omega$ is bounded and $G$ does not satisfy $\Delta_2$ near infinity; if $\Omega$ is bounded and $G$ satisfies $\Delta_2$ near infinity, both sets coincide). See \cite[Remark 3.8]{CFBS}.
\item[(ii)] The domain of $\mathcal J_s$ must be defined through the condition ``$u=0$ a.e.\ outside $\Omega$''. If it is defined instead as the norm closure of $C^\infty_c(\Omega)$ in $W^sL_G(\mathbb{R}^n)$, the limsup inequality may fail for irregular $\Omega$, already for $G(a)=a^2$ \cite[Remark 3.9]{CFBS}.
\item[(iii)] This theorem corrects and extends \cite[Theorem 6.5]{FBS}; see \cite{CFBS}.
\end{enumerate}
\end{rem}

\begin{rem}
If $\Omega$ is bounded, Theorem \ref{teo.comp} gives $W^s_0L_G(\Omega)\subset E_G(\Omega)$ for every $s\in(0,1]$ (for $s=1$, use Lemma \ref{lema.cota}). Hence the effective domains of $\mathcal J_s$ and $\mathcal J$ are contained in $E_G(\Omega)$, which is therefore the natural ambient space.
\end{rem}

\subsection{Minimization problems and the fractional $g$-Laplacian}\label{sec.glap}

Throughout this subsection, $\Omega\subset\mathbb{R}^n$ is a bounded open set (no regularity is assumed), and $f\in E_{\overline{G}}(\Omega)$. No growth condition is imposed on $G$ (neither $\Delta_2$ nor $\nabla_2$). Some of the results below require one of the following conditions:
\begin{enumerate}
\item[(G1)] (regularity) $g=G'$ is continuous, that is, $G\in C^1([0,\infty))$;
\item[(G2)] (strict convexity) $G$ is strictly convex, that is, $g$ is strictly increasing.
\end{enumerate}
Condition (G1) is used only to derive the Euler--Lagrange equation of the fractional problem, and (G2) only for the uniqueness of its minimizers. The limit problem requires neither of them, since $ G_0 $ always satisfies both (Lemma \ref{lema.Gt.reg}). Since $\Omega$ is bounded, $L^\infty(\Omega)\subset E_{\overline{G}}(\Omega)$. If $\overline{G}$ satisfies $\Delta_2$ (equivalently, if $G$ satisfies the $\nabla_2$ condition), then $E_{\overline{G}}(\Omega)=L_{\overline{G}}(\Omega)$; see Remark \ref{rem.f.LG} for $f\in L_{\overline{G}}(\Omega)$ in general.

For $s\in(0,1]$ we consider the functional $\mathcal F_s\colon E_G(\Omega)\to(-\infty,+\infty]$,
\[
\mathcal F_s(u)=\mathcal J_s(u)-\int_\Omega fu\,dx,
\]
where $\mathcal J_1=\mathcal J$. Since $u\mapsto\int_\Omega fu$ is continuous on $E_G(\Omega)$, Theorem \ref{teo.gamma.conv} and \cite[Proposition 6.21]{DalMaso} give
\begin{equation}\label{gamma.F}
\mathcal F_1=\Gamma\text{-}\lim_{s\uparrow1}\mathcal F_s\qquad\text{in }E_G(\Omega).
\end{equation}

\begin{lema}\label{lema.Gt.reg}
For every Young function $G$, $ G_0 \in C^1([0,\infty))$, $ G_0 $ is strictly convex, and $\xi\mapsto G_0 (|\xi|)$ is strictly convex on $\mathbb{R}^n$.
\end{lema}

\begin{proof}
Let $\Lambda(\tau)=\int_{\mathbb S^{n-1}}G(\tau|\omega\cdot e_n|)\,d\mathcal H^{n-1}$. For $a>0$, define $g_0(a)=\Lambda(a)/a$, and set $g_0(0)=0$. As in the proof of Proposition \ref{prop.equiv}, this function is continuous and $ G_0 '=g_0$ on $[0,\infty)$. First, $G(r)/r=\frac1r\int_0^r g$ is strictly increasing on $(0,\infty)$. It is nondecreasing by the monotonicity of $g$. If it were constant, with value $c>0$, on an interval $[r_1,r_2]$ with $0<r_1<r_2$, then $G(r)=cr$ there and hence $g=c$ almost everywhere on $(r_1,r_2)$. Monotonicity gives $g\le c$ almost everywhere on $(0,r_1)$, while $G(r_1)=cr_1$ forces $g=c$ almost everywhere on $(0,r_1)$. This contradicts the right-continuity of $g$ at $0$ and $g(0)=0$. Since
\[
g_0(a)=\int_{\mathbb S^{n-1}}|\omega\cdot e_n|\,\frac{G(a|\omega\cdot e_n|)}{a|\omega\cdot e_n|}\,d\mathcal H^{n-1},\qquad a>0,
\]
(with the integrand equal to $0$ where $\omega\cdot e_n=0$), and $|\omega\cdot e_n|>0$ on a set of positive measure, $g_0$ is strictly increasing, that is, $ G_0 $ is strictly convex. Finally, let $\xi\ne\eta$ in $\mathbb{R}^n$. If $|\xi|\ne|\eta|$, then $ G_0 (|\frac{\xi+\eta}2|)\le G_0 (\frac{|\xi|+|\eta|}2)<\frac12( G_0 (|\xi|)+ G_0 (|\eta|))$. If $|\xi|=|\eta|$, then $|\frac{\xi+\eta}2|<|\xi|$ and the strict inequality follows since $ G_0 $ is strictly increasing.
\end{proof}

\begin{thm}\label{teo.exist}
For every $s\in(0,1]$, $\mathcal F_s$ has a minimizer, and every minimizer $u_s$ satisfies
\[
[u_s]_{s,G}\le\frac{R_0}{\alpha_n},\qquad\mathcal J_s(u_s)\le2A,\qquad\|u_s\|_G\le\frac{C_\Omega R_0}{\alpha_n},
\]
where the constants $R_0$, $A$, $C_\Omega$ (defined in the proof) do not depend on $s$. The minimizer of $\mathcal F_1$ is unique. If $G$ satisfies (G2), the minimizer of $\mathcal F_s$ is unique for every $s\in(0,1)$.
\end{thm}

\begin{proof}
Let $R_\Omega$, $\tilde c_n$ be as in Proposition \ref{prop.poinc.s}, $C_\Omega$ as in Theorem \ref{teo.poinc.unif}, and
\[
A=\int_\Omega\overline{G}\Big(\frac{2R_\Omega}{\tilde c_n}|f|\Big)dx,
\]
which is finite since $f\in E_{\overline{G}}(\Omega)$. Let $u\in W^s_0L_G(\Omega)$ and $\lambda=2/\tilde c_n\ge1$. By Young's inequality, the convexity of $G$ and \eqref{poinc.mod},
\[
\int_\Omega fu\,dx\le\int_\Omega G\Big(\frac{|u|}{\lambda R_\Omega}\Big)dx+\int_\Omega\overline{G}(\lambda R_\Omega|f|)\,dx\le\frac1{\lambda \tilde c_n}\mathcal J_s(u)+A=\frac12\mathcal J_s(u)+A .
\]
Hence
\begin{equation}
\mathcal F_s(u)\ge\frac12\mathcal J_s(u)-A\qquad\text{for every }s\in(0,1]\text{ and }u\in W^s_0L_G(\Omega).
\end{equation}
In particular $\inf\mathcal F_s\in[-A,0]$, and every $u$ with $\mathcal F_s(u)\le0$ satisfies $\mathcal J_s(u)\le2A$. By convexity, this gives $[u]_{s,G}\le R_0:=\max\{1,2A\}$ for $s<1$ and $[u]_{1, G_0 }\le R_0$ for $s=1$. By \eqref{equiv.normas.1}, since $\alpha_n\le1$, and Theorem \ref{teo.poinc.unif},
\begin{equation}\label{cotas.unif}
[u]_{s,G}\le\frac{R_0}{\alpha_n},\qquad\mathcal J_s(u)\le2A,\qquad\|u\|_G\le\frac{C_\Omega R_0}{\alpha_n}\qquad\text{for every }s\in(0,1].
\end{equation}

\emph{Existence for $s<1$.} Let $\{u_k\}$ be a minimizing sequence with $\mathcal F_s(u_k)\le0$. By \eqref{cotas.unif}, $\{u_k\}$ is bounded in $W^sL_G(\mathbb{R}^n)$, since $[u_k]_{s,G}\le R_0/\alpha_n$. By Theorem \ref{teo.comp}, and since the $u_k$ vanish outside $\Omega$, a subsequence converges in $L_G(\Omega)$ and a.e.\ to some $u$. By Fatou's lemma, $\Phi_{G,\nu_s}(D^su)\le\liminf\Phi_{G,\nu_s}(D^su_k)$, so $u\in W^s_0L_G(\Omega)$, and since $\int fu_k\to\int fu$, $u$ is a minimizer.

\emph{Existence for $s=1$.} Let $\{u_k\}$ be a minimizing sequence with $\mathcal F_1(u_k)\le0$, and let $\kappa=\max\{1,2C_\Omega R_0/\alpha_n\}$. Then $\Phi_G(2u_k/\kappa)\le1$ and $\Phi_{ G_0 }(|\nabla u_k|/\kappa)\le2A$. Applying Lemma \ref{lema.cota} to $u_k/\kappa$ gives a uniform bound for $\{u_k/\kappa\}$ in $W^{1/2}L_G(\mathbb{R}^n)$; since $\kappa$ is fixed, the same is true of $\{u_k\}$. By Theorem \ref{teo.comp}, a subsequence converges in $L_G(\Omega)$ to some $u$, and by Lemma \ref{lema.cierre}, $u\in W^1_0L_G(\Omega)$ and $\mathcal J(u)\le\liminf\mathcal J(u_k)$. Hence $u$ is a minimizer.

\emph{Uniqueness.} By Lemma \ref{lema.Gt.reg}, $\xi\mapsto G_0 (|\xi|)$ is strictly convex on $\mathbb{R}^n$. If $G$ satisfies (G2), then $a\mapsto G(|a|)$ is strictly convex on $\mathbb{R}$ (the argument is the same as in Lemma \ref{lema.Gt.reg}). If $u\ne v$ are minimizers of $\mathcal F_s$, then $D^su\ne D^sv$ on a set of positive $\nu$-measure (resp.\ $\nabla u\ne\nabla v$ on a set of positive measure), since otherwise $u-v$ would be constant on $\mathbb{R}^n$, hence zero. Therefore $\mathcal F_s(\frac{u+v}2)<\frac12(\mathcal F_s(u)+\mathcal F_s(v))$, a contradiction.

The bounds follow from \eqref{cotas.unif}, since $\mathcal F_s(u_s)\le\mathcal F_s(0)=0$.
\end{proof}

\begin{thm}\label{teo.conv.min}
For $s\in(0,1)$ let $u_s$ be any minimizer of $\mathcal F_s$, and let $u_1$ be the minimizer of $\mathcal F_1$. Then $u_s\to u_1$ in $E_G(\Omega)$ and $\min\mathcal F_s\to\min\mathcal F_1$ as $s\uparrow1$.
\end{thm}

\begin{proof}
Let $s_k\uparrow1$ and $\kappa=\max\{1,C_\Omega R_0/\alpha_n\}$, with the notation of the previous proof. The functions $v_k=u_{s_k}/\kappa$ satisfy $\Phi_G(v_k)\le1$ and $\Phi_{G,\nu_{s_k}}(D^{s_k}v_k)\le2A$. By Theorem \ref{main2}, a subsequence of $\{u_{s_k}\}$ converges in $L_{G,{\rm loc}}(\mathbb{R}^n)$, hence in $E_G(\Omega)$ since all the functions vanish outside $\Omega$. By \eqref{gamma.F} and the fundamental theorem of $\Gamma$-convergence \cite[Corollary 7.20]{DalMaso}, every limit point minimizes $\mathcal F_1$, and $\min\mathcal F_{s_k}\to\min\mathcal F_1$. Since the minimizer of $\mathcal F_1$ is unique, the limit is $u_1$ and the whole family converges.
\end{proof}

\begin{rem}\label{rem.f.LG}
The assumption $f\in E_{\overline{G}}(\Omega)$ is only used to ensure that $A<\infty$. If $f\in L_{\overline{G}}(\Omega)$ and $\|f\|_{\overline{G}}\le\frac{\tilde c_n}{2R_\Omega}$, then $A\le1$, and Theorems \ref{teo.exist} and \ref{teo.conv.min} hold with the same proofs. No growth condition on $G$ is needed: the superlinear term is provided by the modular Poincar\'e inequality \eqref{poinc.mod}, whose constants do not depend on $s$.
\end{rem}

Next we relate minimizers with solutions of the associated Euler--Lagrange equations. Without the $\Delta_2$ condition the functionals $\mathcal J_s$ are not differentiable on their domains: $\Phi_{G,\nu_s}(D^s(u+\lambda v))$ may be infinite for every $\lambda\ne0$. Following the ideas of \cite{Gossez, FBSp}, we use a variational inequality formulation. For $s\in(0,1)$ let
\[
\mathcal D_s=\{w\in W^s_0L_G(\Omega)\colon\Phi_{G,\nu_s}(D^sw)<\infty\},\qquad\mathcal D_1=\{w\in W^1_0L_G(\Omega)\colon\Phi_{ G_0 }(|\nabla w|)<\infty\},
\]
and, for $\xi\in\mathbb{R}^m$, let $\mathbf g(\xi)=g(|\xi|)\frac{\xi}{|\xi|}$ and $\mathbf g_0(\xi)=g_0(|\xi|)\frac{\xi}{|\xi|}$ ($=0$ if $\xi=0$), where $g_0(a)= G_0 '(a)=\Lambda(a)/a$ is continuous. Under (G1), these are the gradients of the convex $C^1$ functions $G(|\xi|)$ and $ G_0 (|\xi|)$; for $ G_0 $ this holds for every Young function $G$, by Lemma \ref{lema.Gt.reg}.
\medskip

\paragraph{The fractional $g$-Laplacian.} Let $s\in(0,1)$. For $u$ regular enough (for instance, $u\in C^2_c(\mathbb{R}^n)$ and $g$ locally Lipschitz), the \emph{fractional $g$-Laplacian} of $u$ is
\begin{equation}\label{gLap}
(-\Delta_g)^su(x)=2\,\operatorname{p.v.}\int_{\mathbb{R}^n}\mathbf g\big(D^su(x,y)\big)\frac{dy}{|x-y|^{n+s}}=2\lim_{\varepsilon\downarrow0}\int_{\{|x-y|>\varepsilon\}}\mathbf g\big(D^su(x,y)\big)\frac{dy}{|x-y|^{n+s}},
\end{equation}
where $\mathbf g(t)=g(|t|)\frac{t}{|t|}$ ($=0$ if $t=0$) for $t\in\mathbb{R}$. Since $\mathbf g$ is odd and $D^su(y,x)=-D^su(x,y)$, exchanging $x$ and $y$ shows that, for such $u$ and every $v\in C^\infty_c(\mathbb{R}^n)$,
\begin{equation}\label{gLap.weak}
\int_{\mathbb{R}^n}(-\Delta_g)^su\,v\,dx=\iint_{\mathbb{R}^{2n}}\mathbf g(D^su)\,D^sv\,d\nu .
\end{equation}
The right-hand side of \eqref{gLap.weak} makes sense under much weaker assumptions, and we use it as the definition in general: if $u$ is measurable and $g(|D^su|)\in L_{\overline G}(\mathbb{R}^{2n},\nu)$, then $(-\Delta_g)^su$ is the linear functional
\[
\langle(-\Delta_g)^su,v\rangle=\iint_{\mathbb{R}^{2n}}\mathbf g(D^su)\,D^sv\,d\nu,\qquad v\in W^sL_G(\mathbb{R}^n),
\]
which is well defined and continuous on $W^sL_G(\mathbb{R}^n)$ by H\"older's inequality, since $|\langle(-\Delta_g)^su,v\rangle|\le2\|g(|D^su|)\|_{\overline G,\nu}\|D^sv\|_{G,\nu}$. Thus $(-\Delta_g)^su$ is formally the derivative of $u\mapsto\Phi_{G,\nu}(D^su)$, and the factor $1-s$ in \eqref{ec.g} below is the normalization of $\nu_s$, for which $(1-s)(-\Delta_g)^s$ is the natural approximation of $-\operatorname{div}(\mathbf g_0(\nabla\,\cdot))$ as $s\uparrow1$, in the sense of Theorems \ref{main1} and \ref{teo.conv.sol}. When $G(a)=a^p/p$, $(-\Delta_g)^s$ is (up to a factor $2$) the fractional $p$-Laplacian.

\begin{defn}
Let $s\in(0,1)$. We say that $u$ is a \emph{variational solution} of
\begin{equation}\label{ec.g}
(1-s)(-\Delta_g)^su=f\quad\text{in }\Omega,\qquad u=0\quad\text{in }\mathbb{R}^n\setminus\Omega,
\end{equation}
if $u\in\mathcal D_s$, $g(|D^su|)\in L_{\overline{G}}(\mathbb{R}^{2n},\nu)$ and
\begin{equation}\label{VI}
(1-s)\iint_{\mathbb{R}^{2n}}\mathbf g(D^su)\,(D^sw-D^su)\,d\nu\ge\int_\Omega f(w-u)\,dx\qquad\text{for every }w\in\mathcal D_s .
\end{equation}
Analogously, $u$ is a variational solution of
\begin{equation}\label{ec.local}
-\operatorname{div}\big(\mathbf g_0(\nabla u)\big)=f\quad\text{in }\Omega,\qquad u=0\quad\text{on }\partial\Omega,
\end{equation}
if $u\in\mathcal D_1$, $g_0(|\nabla u|)\in L_{\overline{ G_0 }}(\mathbb{R}^n)$ and $\int\mathbf g_0(\nabla u)\cdot(\nabla w-\nabla u)\,dx\ge\int_\Omega f(w-u)\,dx$ for every $w\in\mathcal D_1$.
\end{defn}

By H\"older's inequality, the integrals in \eqref{VI} are finite. By \eqref{gLap.weak}, the left-hand side of \eqref{VI} is $(1-s)\langle(-\Delta_g)^su,w-u\rangle$. The operator \eqref{gLap} coincides with the fractional $g$-Laplacian of \cite{FBS, FBSp} up to a positive normalization constant.

\begin{thm}\label{teo.equiv}
Let $s\in(0,1]$, and if $s<1$ assume (G1). Let $u\in\mathcal D_s$. Then $u$ minimizes $\mathcal F_s$ if and only if $u$ is a variational solution of \eqref{ec.g} (resp.\ of \eqref{ec.local} if $s=1$). Moreover, in that case $u$ is a weak solution in the following sense: for every $v\in W^s_0E_G(\Omega)$, and in particular for every $v\in C^\infty_c(\Omega)$,
\[
(1-s)\iint_{\mathbb{R}^{2n}}\mathbf g(D^su)\,D^sv\,d\nu=\int_\Omega fv\,dx\qquad\Big(\text{resp. }\int\mathbf g_0(\nabla u)\cdot\nabla v\,dx=\int_\Omega fv\,dx\Big).
\]
\end{thm}

\begin{proof}
We give the proof for $s<1$; the case $s=1$ is identical, replacing $G$, $D^s$, $\nu$ and the factor $1-s$ by $ G_0 $, $\nabla$, the Lebesgue measure and $1$. Let $\phi(a)=G(|a|)$, so that $\phi'=\mathbf g$ by (G1). For $s=1$, $ G_0 \in C^1$ by Lemma \ref{lema.Gt.reg}, so no assumption is needed.

\emph{Minimizer $\Rightarrow$ variational solution.} Let $w\in\mathcal D_s$. By convexity, $u+\theta(w-u)=(1-\theta)u+\theta w\in\mathcal D_s$ for $\theta\in(0,1]$, and by minimality
\[
(1-s)\iint q_\theta\,d\nu\ge\int_\Omega f(w-u)\,dx,\qquad q_\theta=\frac{\phi(D^su+\theta(D^sw-D^su))-\phi(D^su)}{\theta}.
\]
Since $\phi$ is convex, $q_\theta$ decreases as $\theta\downarrow0$ to $\mathbf g(D^su)(D^sw-D^su)$, and $q_\theta\le q_1=\phi(D^sw)-\phi(D^su)\in L^1(\nu)$. By the monotone convergence theorem,
\begin{equation}\label{VI.debil}
(1-s)\iint\mathbf g(D^su)(D^sw-D^su)\,d\nu\ge\int_\Omega f(w-u)\,dx\qquad\text{for every }w\in\mathcal D_s,
\end{equation}
where the left-hand side is a priori in $[-\infty,\infty)$. Taking $w=0$ in \eqref{VI.debil},
\[
(1-s)\iint g(|D^su|)|D^su|\,d\nu\le\int_\Omega fu\,dx<\infty .
\]
Since $\overline{G}(g(a))=ag(a)-G(a)\le ag(a)$, we get $\iint \overline{G}(g(|D^su|))\,d\nu<\infty$, so $g(|D^su|)\in L_{\overline{G}}(\nu)$, and \eqref{VI.debil} is \eqref{VI}.

\emph{Variational solution $\Rightarrow$ minimizer.} By convexity, $\phi(b)\ge\phi(a)+\mathbf g(a)(b-a)$. Hence, for $w\in\mathcal D_s$,
\[
\mathcal F_s(w)-\mathcal F_s(u)\ge(1-s)\iint\mathbf g(D^su)(D^sw-D^su)\,d\nu-\int_\Omega f(w-u)\,dx\ge0,
\]
and $\mathcal F_s(w)=+\infty$ if $w\notin\mathcal D_s$.

\emph{Weak form.} Let $\ell(z)=(1-s)\iint\mathbf g(D^su)D^sz\,d\nu-\int_\Omega fz\,dx$, which is a linear functional on $W^s_0L_G(\Omega)$ by H\"older's inequality. By \eqref{VI}, $\ell(w)\ge\ell(u)$ for every $w\in\mathcal D_s$. If $v\in W^s_0E_G(\Omega)$, then $\lambda v\in\mathcal D_s$ for every $\lambda\in\mathbb{R}$, so $\lambda\ell(v)\ge\ell(u)$ for every $\lambda\in\mathbb{R}$, which forces $\ell(v)=0$. Finally, $C^\infty_c(\Omega)\subset W^s_0E_G(\Omega)$ by Lemma \ref{lema.cota}.
\end{proof}

\begin{rem}
If $G$ satisfies $\Delta_2$, then $\mathcal D_s=W^s_0L_G(\Omega)=W^s_0E_G(\Omega)$ is a linear space, and variational solutions coincide with weak solutions in the usual sense. In general, the weak form in Theorem \ref{teo.equiv} only allows test functions in $W^s_0E_G(\Omega)$, and the variational inequality \eqref{VI} is the appropriate substitute.
\end{rem}

Combining Theorems \ref{teo.exist}, \ref{teo.conv.min} and \ref{teo.equiv} we obtain the main result of this subsection, which extends the main result of \cite[Section 6.3]{FBS} to Young functions not satisfying $\Delta_2$.

\begin{thm}\label{teo.conv.sol}
There exists a unique variational solution $u$ of \eqref{ec.local}. Assume (G1). Then, for every $s\in(0,1)$, there exists a variational solution $u_s$ of \eqref{ec.g}, which is unique if moreover $G$ satisfies (G2). For any choice of the $u_s$, $u_s\to u$ in $E_G(\Omega)$ as $s\uparrow1$, and $u$ satisfies
\[
\int_\Omega\mathbf g_0(\nabla u)\cdot\nabla v\,dx=\int_\Omega fv\,dx\qquad\text{for every }v\in C^\infty_c(\Omega),
\]
that is, $u$ is a distributional solution of \eqref{ec.local}.
\end{thm}

\begin{proof}
By Theorem \ref{teo.equiv}, variational solutions are exactly the minimizers of $\mathcal F_s$, which exist by Theorem \ref{teo.exist}, and are unique for $s=1$, and for $s<1$ under (G2). The convergence follows from Theorem \ref{teo.conv.min}, and the last assertion from Theorem \ref{teo.equiv} with $s=1$.
\end{proof}

\end{document}